# An Operational Scheduling Framework for Tanker-based Water Distribution System under Uncertainty

Abhilasha Maheshwari[1#], Shamik Misra[2#], Ravindra Gudi[3*], Senthilmurugan Subbiah[4], Chrysi Laspidou[5]

[1]Department of Chemical Engineering, Indian Institute of Technology Jodhpur, Jodhpur, Rajasthan - 320347, India

[2]Department of Chemical Engineering, Indian Institute of Technology Tirupati, Tirupati, Andhra Pradesh, India

[3] Department of Chemical Engineering, Indian Institute of Technology Bombay, Powai, Mumbai - 400076, India

[4] Department of Chemical Engineering, Indian Institute of Technology Guwahati, Assam 781039, India

[5] Department of Civil Engineering, University of Thessaly, Pedion Areos, 38334 Volos, Greece

*To whom correspondence should be addressed: Prof. Ravindra D. Gudi, Department of Chemical Engineering, CAD Centre, IIT Bombay, Powai, Mumbai, India – 400076, (phone: +91 (22) 2576 7231, e-mail: ravigudi@iitb.ac.in)

[#] Authors have contributed equally to this work.

**Abstract**

Tanker water systems play critical role in providing adequate service to meet potable water demands in the face of acute water crisis in many cities globally. Managing tanker movements among the supply and demand sides requires an efficient scheduling framework that could promote economic feasibility, ensure timely delivery, and avoid water wastage. However, to realize such a sustainable water supply operation, inherent uncertainties related to consumer demand and tanker travel time need to accounted in the operational scheduling. Herein, a two-stage stochastic optimization model with a recourse approach is developed for scheduling and optimization of tanker based water supply and treatment facility operations under uncertainty. The uncertain water demands and tanker travel times are combinedly modelled in a computationally efficient manner using a hybrid Monte Carlo simulation and scenario tree approach. The maximum demand fulfillment, limited extraction of groundwater, and timely delivery of quality water are enforced through a set of constraints to achieve sustainable operation. A representative urban case study is demonstrated, results are discussed for two uncertainty cases (i) only demand, and (ii) integrated demand-travel time. Value of stochastic solution over expected value and perfect information model solutions are analyzed and features of the framework for informed decision-making are discussed.

1. Introduction

Providing affordable access to clean and pure water to all is one of the main objectives of Sustainable Development Goals (SDG 6) [1]. Provision of such a service is mainly via two main methods, namely, piped water distribution network and tanker-based water supply systems. To this end, challenges related to water supply in urban areas include (i) low reliability of piped distribution systems supply[2, 3], and (ii) lack of full coverage by piped water distribution network in peri-urban areas[4]. Piped water supply being spatially and temporally uneven across the city, tanker-based water supply is one of the most prominent forms of alternative arrangements for water distribution in cities of various countries like India, Italy, Malaysia, Kenya etc. [5-7]. These tanker water distribution systems facilitate residents in peri-urban areas who are not connected to piped network infrastructure and also in urban areas where people are facing water shortages due to highly intermittent nature of piped distribution systems. The main advantage is that the tanker water distribution system can be quickly implemented to cater to a large number of customers with variable demands and water quality requirements with significantly less infrastructural challenges as compared to piped distribution network.

Statistics reports that only 49% of households in India have access to a drinking water tap from a treated water source in current scenario [8]. Such conditions have increased the dependence of public on tanker-based distribution systems. Not only limited to the developing countries, this tanker water supply system is also prevalent in developed nations like Canada, where 13% of households use tanker water as a primary potable water supply source [9]. Therefore, water delivery by tankers is increasingly being considered by local government bodies in urban water management policies and related decision-making for city water supply requirements, e.g. Delhi and Chennai [10, 11]. Consequently, tanker water supply systems are expected to play an important role and receive a growing attention in urban water governance [12].

Acknowledging the water supply through tankers in urban areas, significant attention has been drawn by few researchers through literature studies in recent years. Most of these studies have described the inadequate services in prevalent tanker water supply systems [13, 14] and the economic evaluation of tanker water supply considering the policy and environmental factors [15]. Many of these studies report challenges in relation to the poor water quality, informal settlements of vendors, inadequate service features and high price of water charged to consumers for delivery [16, 17].

Furthermore, [6] conducted a detailed survey of tanker water market in the city of Jordan and developed a simulation model to estimate tanker water consumption behavior (demand pattern) of commercial consumers. They also analyzed the influence of spatial factors such as water source location and pipe network supply duration on commercial demand fulfillment. Nevertheless, research gaps pertaining to systematic operation of tanker water supply e.g. efficient management of large fleet of tanker trucks, proper coordination among spatially distributed water sources and customers spanning large regions, remains largely anecdotal in the literature [14].

In this direction, (Maheshwari et al. 2020) [18] developed a short-term planning formulation for optimal operation of tanker based water distribution systems. The study presented a MILP optimization formulation minimizing the total operating cost and estimating the optimal water treatment and tanker water distribution plan considering various constraints related to (i) treatment facilities operation, (ii) different consumer type and their demand fulfillment for different water quality products, and (iii) tanker availability. The tanker movements schedule and source-consumer-tanker association matrices generated from the framework promoted efficient coordination between tanker suppliers, water treatment plants (WTPs) operation, and consumers for timely delivery of quality water as well as optimal use of available water resources. The formulation was though based on the assumed known water demands from each consumer and fixed travel time calculated from distance to be travelled.

However, tanker water supply management in urban areas is also faced with multiple type of uncertainties such as short-term variation in water demands, time delays in delivery due to traffic congestions, seasonal water availability and long-term uncertainty in climate conditions, population growth, urbanization etc. [5]. Thus, the challenge of establishing a balance of the trade-off between water supply shortage risks and transportation cost, in the entire supply system is compounded by multiple critical uncertainties. This makes the deterministic scheduling solution invalid in the real-time, when these variables exhibit deviations from the expected values. Therefore, these aspects pose great significance and inevitably needs to be considered in optimal scheduling of tanker-based water distribution system operations under uncertain situations to ensure efficient and reliable performance. Nonetheless, no existing literature deals with the optimal operational scheduling and management of the tanker water supply systems under various uncertainty factors, such as tanker travel time and water demand. A rudimentary version accounting only for the demand uncertainty in the tanker water supply problem has been presented with an abstract formulation in a recent work (Maheshwari et al. 2022) [19]. The work reported a very small example study, with only focus only on demonstrating the value of the stochastic solution over expected value solution. However, various salient aspects of this problem such as optimal use of water sources, tanker delivery scheduling, maximum tanker capacity utilization, source-consumer mapping etc. are included and additionally analyzed in this paper. Furthermore, (Maheshwari et al. 2022) [19] neither discussed nor had the scope to delineate on the complexities arising in formulation due to integration of demand and travel time uncertainty in the tanker system operations, necessitating a more meticulous depiction. . A complete and realistic problem representation along with a detailed and computationally efficient formulation is thus discussed in this paper for addressing integrated uncertainty problem.

This study incorporates several important nuances of tanker-based water distribution system, such as (i) different throughput scale of WTPs, (ii) raw and treated water reservoir capacity limitations and targets, (iii) quality of water for different purposes, (iv) variety of tanker distribution capacity, (v) consumer water consumption pattern, (vi) delays in travel time of tankers due to traffic, roadblocks etc. in the proposed optimization formulation. Addressing these peculiarities and uncertainties herein, we propose a stochastic mixed integer linear programming (MILP) model with hybrid approach of scenario generation for the optimal scheduling of a tanker supply system considering various operational constraints. In this paper, the uncertainties in tanker travel time and water demands are integrated in a hybrid manner, coupling the Monte Carlo simulation (MCS) and scenario tree (ST) methods of stochastic modelling. Thus, the main contributions of this work are:

(1) Development of a novel operational scheduling framework that incorporates the uncertainties inherent to tanker water distribution system and provides optimal transportation schedule to minimize the total operating cost.

(2) The proposed two-stage stochastic MILP model includes several operational, logistic and ground water extraction constraints of tanker water supply system, including different treatment capacities of WTPs, raw and treated water product storage reservoir limitations, source and consumer association selection for different products, various tanker distribution capacities, etc. The decisions regarding (i) water treatment plant operations, (ii) raw (untreated) water supply from different sources to WTPs, and (iii) delivery scheduling of treated water to different consumers in different zones of the city area, are solved concurrently.

(3) The sample average approximation technique (SAA) coupled with the Monte Carlo method is used to make the problem computationally tractable, reducing the stochastic model into an equivalent deterministic optimization problem with multiple scenarios.

In this paper, we demonstrate the application of proposed framework and solution approach on a representative tanker water supply system case study, typical to Indian cities and the effects of demand and transportation time uncertainty on the optimal scheduling of tanker-based water distribution system is discussed in detail.

Section 2 elucidates the tanker-based water distribution system problem and defines the objective, decision variables, input parameters, uncertainty, and other assumptions of this stochastic optimization problem. In Section 3, a brief discussion about the two-stage stochastic programming approach is presented and the abstract form of the mathematical optimization problem formulation is described. Section 4 presents a representative case study to demonstrate the potential and features of the proposed framework and discusses several important aspects of operational scheduling under uncertainty in tanker-based water distribution systems in semi/urban regions. A detailed comparison of optimization results for stochastic formulation with the expected value solution for water demands is also analyzed and illustrated to understand the influence of the uncertainty on various cost components. Section 5 provides the main conclusions along with practical implications of this work in safe water supply.

## 2. Problem Overview

This section delineates the tanker water supply operations scheduling problem under uncertainty and proposed hybrid model for uncertainty modelling.

**2.1 Problem definition**

The tanker water supply system typically includes several water sources located in different areas, advanced water treatment plants (WTPs) and customers. The water sources include ground water borings, rivers, ponds, and lakes, from where raw water is supplied in the tankers to water treatment plants and consumers after an appropriate treatment for the intended purpose. The water for domestic purposes such as drinking, cooking, bathing, and cleaning etc. from fresh water sources (e.g., rivers, lakes etc.) requires only chlorine-based disinfection in the tanker truck after primary treatment at the source site itself [20]. On the other hand, water pumped from ground water sources is first treated in water treatment plants so as to treat to the required quality levels [17]. The raw water is treated using several membranes based advanced wastewater treatment technologies in WTPs (e.g. ultra/nano filtration, reverse osmosis etc.) to reduce high levels of suspended solids (TSS) and dissolved salts (TDS) in raw water [21]. In this line, herein three type of water quality states are considered: (i) raw (untreated) water (RW) extracted from groundwater bore hole sources and sent to WTPs, (ii) domestic purpose water (DPW) after chlorine-based disinfection on tanker trucks, and (iii) ultra-pure drinking water (UPDW) through advanced treatment at WTPs. However, as formulated, the optimization model is sufficiently flexible to incorporate an arbitrary water quality state, as required for the system.

Furthermore, as can be seen in Figure 1, there are four type of tanker movements happening in the entire system, (i) for supplying RW from groundwater sites to WTPs for suitable advance water treatment, (ii) for supplying DPW directly from fresh-water source sites to customers post disinfection treatment for required time period on tanker trucks itself, (iii) for delivering treated water products from WTPs to customers, and (iv) empty tankers returning to nearest region tanker depot to start new trip. In addition, the system also considers different tanker type (in terms of volume and material of construction) to supply water in each section. The optimal selection of tanker type in the system is also constrained to geographical location of consumers, demands (both water quality and quantity), availability and transportation costs in each section. Also, uncertainty in consumer water demands (often forecasted based on historical usage pattern), relatively inadequate road infrastructure and traffic conditions etc. significantly impact the tanker movements and system operations for demand fulfillments. Further, as we know that for the stochastic optimization problem, when proper routing decisions are incorporated, this problem becomes intractable. So, we have proposed to first solve a high-level problem to get decisions about production targets to the treatment facility as well as the sourcing targets to the freshwater and groundwater sources, and a better customer-source-vehicle suitability. Once these targets are fixed from this high-level scheduling, a proper vehicle routing with details of congestions in different route segments can be done at the second level scheduling. However, this later aspect is not a subject matter of the current paper.

Therefore, the main aim of this paper is to address the below challenges specific to stochastic optimization in operational scheduling of tanker-based water distribution system in urban (peri-urban) areas:

(i) Rigorous incorporation of uncertain behavior of water demands in the formulation so as to obtain a computationally tractable optimization problem

(ii) incorporation of uncertainty in tanker travel time and delays while scheduling tanker movements for water supply to the consumers

(iii) ensuring appropriate treatment of raw water (depending on type of water source) and quality of water supplied to the consumers for intended purposes

(iv) Minimize transportation costs in the system operations (short-term planning) by optimally mapping tankers (based on availability and suitability) to consumers for various water demands and volumes

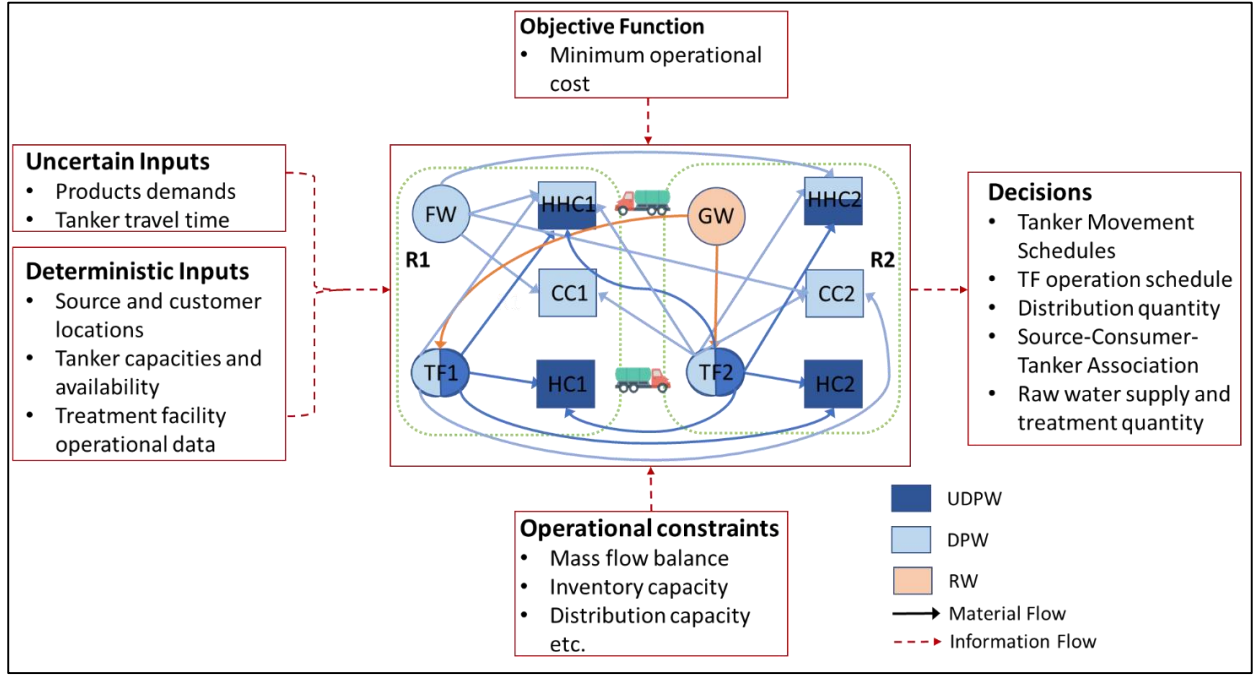

Figure 1: Schematic representation of the stochastic problem

Thus, we describe the problem as follows: the entire city area is split into a number of water supply regions ($r \in R$) depending on the geographical locations of available water sources, WTPs, and consumers. Given this location, a set ($s \in S$) of water sources in a region is used to supply different quality of water ($p \in P$) either to a set of consumers ($c \in C$), or to a set of water treatment plant ($s' \in WTP \subseteq S$) in or across the regions. Furthermore, the volume of water pumped from groundwater sources is constrained by the maximum permissible extraction limit ($S_{max}$) usually specified by the central groundwater board to allow ground water table reclamation. Accounting for water treatment in treatment plants where inventories ($i \in I$) of untreated raw water and treated water are maintained in separate reservoirs. The volume of water in these inventories is limited by the fixed physical capacity of the reservoir tanks. A set of tanker vehicles ($v \in V$) is considered to transport water from source sites to consumer locations and treatment plants. Each region is assumed to have a fixed number of tankers availability and each tanker has a specific capacity, and any extra tankers to be hired have a fixed hiring cost for each type. The operating costs for tanker movements, transporting water from sources to consumers is assumed to be directly proportional to the distance to be traveled and tanker capacities. A distribution time slot is considered in the schedule for each consumer after tankers reaching the destination location of consumer groups in the region.

For efficient management of water resources and timely delivery of tankers with quality water, system managers need to make the tanker water delivery plan according to the supplied quantity of water from sources, treatment facilities, forecasted consumer demands, and the number of available tankers for supply. The water delivery plan includes the volume of raw water to be supplied from each groundwater source to WTPs and treated water to consumers, the required capacity and type of tankers, delivery timings, and other information. Also, since the consumer demands are uncertain, any shortage in the water supply to consumers will result in a penalty cost. Similarly, an extra penalty cost is incurred for any surplus supply of water leading to water wastage. Furthermore, a tanker can make more than one round trip per day; however, traffic congestions leading to uncertainty in transit travel time can reduce the tanker availability in each region.

Hence, considering the uncertainty in the demands and travel time, the main objective is to determine optimal distribution plan and treatment facility operational schedule that minimizes total water transportation cost while ensuring minimum demand shortfall. Thus, the characteristics considered herein for the overall optimization problem are as follows:

**Known Inputs (parameters)**:

- Water sources site and treatment plant locations
- Tanker vehicles: transportation cost parameters, average speed, initial number availability in each region
- Distance between water sources, WTPs and consumers
- WTPs: reservoir capacities, plant throughput capacity, water recovery rate of treatment process
- Demand shortfall/surplus penalty cost parameters

**Uncertain Parameters**:

- Consumer water demands (for each water quality product)
- Tanker travel time from different sources to consumers and WTPs

**Decisions**:

- Raw water and treated water delivery plan from sources to WTPs and consumers respectively. The plan includes both quantity of water and tanker type
- WTP operation schedule for entire scheduling horizon
- Tanker movement schedule for entire scheduling horizon
- Number and type of extra tankers that could be leased at the beginning of scheduling horizon to respond to effects of uncertainty

The schematic representation in Figure 1 summarizes the framework of the optimization model with the considered input parameters, the objective function, the constraints, and the outputs for decision support in operational scheduling.

## 2.2 Stochastic Optimization

The fundamental concept of two-stage stochastic programming is recourse, which allows compensative actions to be taken in the second-stage decisions based on the impacts of uncertainty surfaced after first stage decisions [22][23]. Further nuances of the stochastic optimization problem are described in Hu and Hu (2018), Grossmann (2021)[24,25]. Many researchers have further extended these studies by incorporating the Sample Average Approximation (SAA) [26] method to address the computationally challenging task of infinite possible scenarios of uncertainty realization in stochastic optimization problems [27]. Furthermore, the scenario tree (ST) and Monte Carlo simulation (MCS) are two main modelling techniques that have been used in previous studies to incorporate uncertainty in a stochastic optimization model [28][29][30]. The former technique approximates uncertainty using a finite number of discrete scenarios with the corresponding probability of realization. On the other hand, the latter incorporates uncertainty by generating a sufficiently large number of random scenarios from a continuous probability distribution function (determined from historical time series data), each scenario being equally likely. Accordingly, ST method is computationally efficient only for convex problems and with small number of decision variables[26].

In this paper, MCS method is used to represent uncertainty in water demands (probability density functions is determined based on available historic water consumption profile), leading to a more realistic representation of the demand uncertainties. On the other hand, uncertainty in tankers availability (distribution capacity) due to delays in travel time while supplying water from sources to WTPs or WTPs to consumers can be represented by possible discrete probabilistic scenarios using scenario tree approach, without increasing the computational complexity of the proposed model. Thus, here we have employed a hybrid approach for uncertainty modelling in the optimization formulation, to perform the operational scheduling of tanker-based water distribution systems. Therefore, we present results for both cases, (i) only demand uncertainty using MCS approach, and (ii) demand uncertainty integrated with travel time uncertainty using hybrid modelling approach, to present practical insights on handling the supply system operations under different sets of constraints. The algorithm structure for uncertainty modelling in both cases is shown in Fig.2(a).

*Remark: There are two kinds of problems that are encountered in the tanker-based water distribution problem, namely, inventory routing and vehicle routing. While the vehicle routing problem has features that are more of a continuum nature and needs a good level of granularity in the decision making; the inventory routing problem which is the focus of this paper, has relatively large temporal and spatial scales which increases the complexity, and compels the need for a relatively abstracted representation, to make the problem tractable and computationally efficient. Thus, the formulation needs to consider the following two aspects: (i) ease of interpretation and translation of the decision making into appropriate advisories for planning of the water supply through tankers, and (ii) aspects related to solving and arriving at decisions by decomposing the complexity of the problem. Thus, the main reason for considering only peak and non-peak time periods for the travel time (high and nominal) as two scenarios is primarily to ease the interpretation of the results of this complex decision-making with a tractable formulation where uncertainty can be accommodated; it must be mentioned that the approach is generic enough to be extended to more than two scenarios.*

*Furthermore, in the inventory routing problem, the main decisions that we want to get from this framework is customer, source, and vehicle (tankers) suitability. However, from the complete implementation perspective, the overall decision making needs to be decomposed into two levels (namely, inventory and vehicle routing), where the inventory routing problem will be solved at higher level to first get this suitability and the in the next step the vehicle routing problem, which*

*will take account of specific route segments and detailed congestion periods will be solved [31-34]. Hence, to decide on the required number of tanker vehicles and suitability between customer and sources, the higher-level approximation would be more efficient to solve to optimality without increasing the computational complexity. Further, lower-level problem consists of an efficient detailed vehicle routing problem which can accept customer-source-vehicle targets from the upper level, based on which the routes to reach each customer could be designed, and the vehicle allocation could be done. So, we have proposed to first solve a high-level problem to get decisions about production targets to the treatment facility as well as the sourcing targets to the freshwater and groundwater sources, and a better customer-source-vehicle suitability. Once these targets are fixed from this high-level planning, a detailed vehicle routing with details of congestions in different route segments can be done at the second level scheduling. This bi-level scheme of optimal planning under uncertainty is depicted in the Fig. 2(b). The dark-colored portion of the Fig. 2 (b) represents the subject matter of this manuscript, and the grayed part represents the ongoing research.*

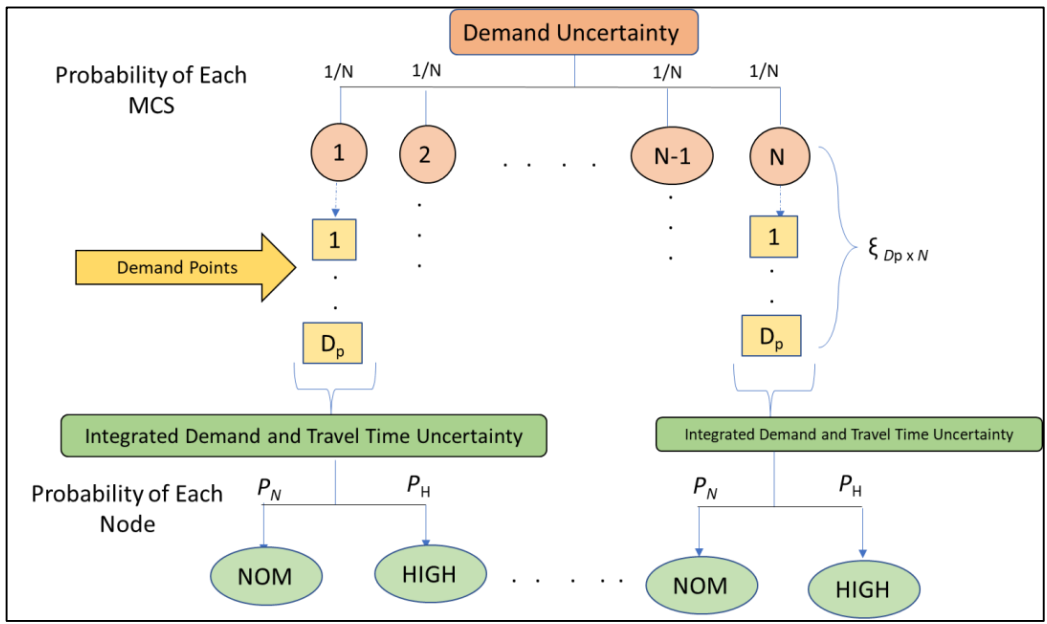

(a)

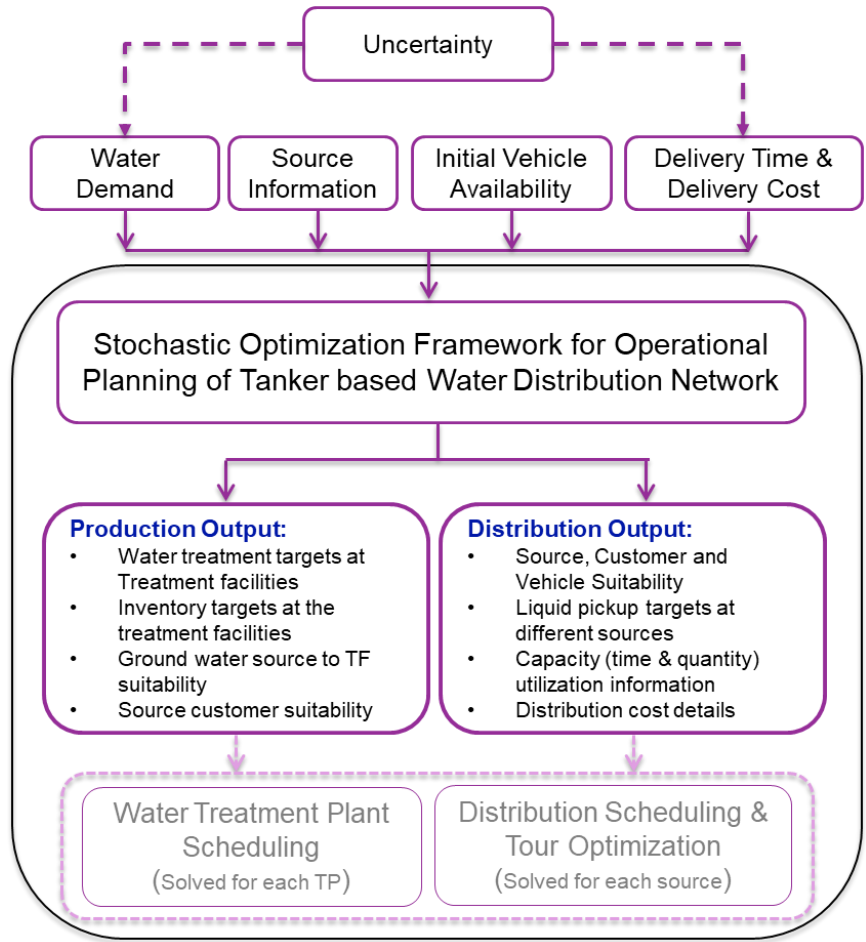

(b)

Fig.2: (a) Algorithm structure for hybrid scenario generation, (b) proposed bilevel framework to achieve the optimal tanker-based water distribution scheme

### 3. Optimization Formulation

A two-stage stochastic MILP formulation with recourse action is developed in this section for addressing the operational scheduling of tanker-based water distribution systems under water demand and tanker travel time uncertainty. The deterministic model for tanker water supply scheduling is presented in (Maheshwari et al. 2020) [18]. Here, a further extension of the model to stochastic case, a detailed mathematical formulation of the equivalent MILP optimization problem is described next, and nomenclature list is provided in the appendix with this paper.

In the proposed two-stage stochastic optimization model, we assumed the decisions of volume of water transported from ground water sources to WTPs, quantity of water contributed by each source to fulfil consumers demand, and extra tanker distribution capacity (initial total number of tankers to start with) as first-stage decision variables, which are to be made before the realization of uncertainty. Decisions related to WTP operation variables and demand shortfall/surplus variables, are considered as second-stage variables. It is to be remarked that while the intent here is above higher level decisions, there are also opportunities for recourse action possible at operator level to take care of some uncertainty in real-time.

### 3.1 Mathematical Formulation:

The main objective of above-described optimization problem is to minimize the total expected operating cost of the tanker water supply system considering various operational constraints over the scheduling horizon. The first stage model is thus formulated as:

$$\underset{x_I}{Min} \; f(x_I) + E(Q(x_I, \xi))$$

s.t.
$$h_I(x_I) = 0$$
$$g_I(x_I) \leq 0 \qquad (1)$$

In this first stage, the model determines the volume of water to be supplied by each water sources to consumer and treatment facilities and the capacity of corresponding tankers required to carry out the distribution task. The first and second terms in the objective function in Eq.1 denote the operating cost corresponding to first-stage decisions, and expected operating cost corresponding to second-stage (recourse) decisions, (given by summation of the demand shortfall/penalty cost in each scenario with its known probability ($\pi_k$)) respectively. The objective function as described in Eq.1 minimizes the expected total operating cost. This cost includes the transportation costs of tankers for water supply from each source to consumers, from GW sources to WTPs, and penalties for violation of buffer and target capacities in raw and treated water reservoirs respectively for smooth operation in the treatment facilities. All the costs in objective function are estimated by linear relationships of the decision variables with corresponding cost parameters. The variables group ($x_I$) in the first term of the objective function are independent of the realization of uncertainty and $\xi$ represents the vector of uncertain parameters in the problem.

Furthermore, having chosen the values of the variables group ($x_I$) in the first stage decision, the future expected cost with the recourse action when the uncertain event takes place is decided in second stage of the model.

The second-stage recourse function is formulated as:

$$Q(x_{II}, y_{II}, \xi) = \min \left\{ \sum_c^C \sum_{p \in P^F}^P \sum_{t=1}^{NT} \left( q_p^+ \Delta dem_{c,p,t}^+(\xi) + q_p^- \Delta dem_{c,p,t}^-(\xi) \right) \right\}$$

s.t.

$$h_{II}(x_I, x_{II}, y_{II}, \xi) = 0$$

$$g_{II}(x_I, x_{II}, y_{II}, \xi) \leq 0 \qquad (2)$$

For a given scenario $\xi_k$, the consumer demand is known in the second stage. Therefore, based on the supplied water quantity decision obtained in the first stage, both the demand shortage and the surplus amount of each consumer and product can be determined. This forms the second-stage total costs as shown in Eq. 2, which is formulated by summing the penalties for failure/surplus in meeting the demands of all consumers for all products in all demand periods. Thus, in the second stage model, for a given scenario $\xi_k$, the recourse function minimizes this cost based on the first stage solution $x_I$. The subscript $k$ in second stage variables indicates decisions under different scenarios of uncertainty realization. Both the shortage and the surplus quantity of each product and consumer are assumed to be non-negative variables at all time periods.

### 3.2 Equivalent Deterministic MILP Model

The above described two-stage stochastic recourse problem can also be formulated as an equivalent deterministic MILP optimization model if the realizations of uncertainty scenario set $\xi$ can be managed to a finite number. To this end, a finite number of scenarios (*N*) are generated using the random sampling method based on Monte Carlo simulations, as described earlier in section 2.2. All the *N* realizations of the uncertain demand parameter in the random vector $\xi$ from this sampling method are based on the probability distribution (from historical time series behaviour) and are independent and identically distributed (*iid*) random variables with equal probabilities of *1/N*. The number of scenarios (*N*) have to determined suitably through repeated computational tests until objective function value converges within the reasonable prescribed limits. Thus, assuming sufficiently large samples, the expected second-stage cost can be then approximated using SAA technique and equivalent problem can be written as follows:

$$\text{Min } f(x_I) + \frac{1}{N}\sum_{k=1}^{K} Q(x_{II}, y_{II}, \xi_k)$$

$$h_I(x_I) = 0$$

$$g_I(x_I) \leq 0$$

$$h_{II}(x_I, x_{II}, y_{II}, \xi_k) = 0 \quad \forall\, k \in K$$

$$g_{II}(x_I, x_{II}, y_{II}, \xi_k) \leq 0 \quad \forall\, k \in K \tag{3}$$

Accordingly, the recourse function can now be re-formulated as following using discrete scenarios:

$$\text{Min } \sum_{k}^{K} \pi_k \sum_{c}^{C} \sum_{p \in P^F}^{P} \sum_{t=1}^{NT} (q_p^+ \Delta dem_{c,p,t,k}^+ + q_p^- \Delta dem_{c,p,t,k}^-)$$

(4)

Thus, in the above equivalent model, the second-stage (operational) constraints have to be satisfied for each scenario, where the optimal value of operational variables can be different for each scenario. However, the first-stage constraints will be satisfied always for all scenarios and the solved water supply plan remains applicable in any possible conditions represented by scenarios.

A detailed mathematical formulation of equivalent deterministic problem of the two-stage stochastic problem with recourse approach, starting with the constraints related to treatment facilities operations, consumer demand fulfillments, tanker capacity utilization, uncertainty modelling followed by the objective function of the problem as the total expected operating cost is presented next. It is to be noted here that all the variables in the following MILP formulation are considered to be positive unless specified otherwise and colon symbol (:) is used to represent "such that", as typically used in mathematical formulations. To develop a sufficiently rigorous model to assist in decision making at the planning layer and generate realistic targets for scheduling water tanker movements, a uniform and hourly discretization of the planning horizon is adopted in the following MILP optimization framework.

*Model Assumptions:*

(i) Each tanker serves only for one consumer group at a time. In other words, a tanker distributes entire water in one location only and not allowed to be partially distributing in different locations.

(ii) It is assumed that the demands for domestic purpose water (DPW) can be fulfilled either by on-tanker disinfection of water from any fresh water sources or by suitable treatment of water from ground water sources at TF. However, demands for ultra-pure drinking water (UPDW) can be fulfilled only by TF after advanced water treatment.

(iii) This paper considers the influence of traffic congestions in urban areas on the travel time of tankers. For simplicity, we assumed that there are only two levels of travel time, namely the nominal time (based on distance and average speed of tankers) and over time of fixed duration (2h).

**TF Operations Constraints**

The equations (5-9) combinedly indicates the water treatment plant operation constraints including switching to shut down mode (in case no treatment is required) or to continue operation along with shut down for a minimum period in each mode after switching, to be economically effective in the treatment process. The binary variables 0-1 are used judiciously in combination with continuous variables in these constraints to account for transition of treatment operations (shut down/restart periods) in TF scheduling.

$xSUp_{s,t,k} - xSDn_{s,t,k} = yOp_{s,t,k} - yOp_{s,t-1,k}$

$\forall s(:STy_s = 'TF'), t(:t = 2..NT), k$ (5)

$xSUp_{s,t,k} - xSDn_{s,t,k} = yOp_{s,t,k} - Op_{s,t-1}^{ini}$

$\forall s(:STy_s = 'TF'), t(:t = 1), k$ (6)

$\sum_{m=t-T_s^{UT}+1:(t-T_s^{UT}+1)>0}^{t} xSUp_{s,m,k} \leq yOp_{s,t,k}$

$\forall s(:STy_s = 'TF'), t(:t = 1..NT), k$ (7)

$$\sum_{m=t-T_s^{DT}+1:(t-T_s^{DT}+1)>0}^{t} xSDn_{s,m,k} \leq 1 - yOp_{s,t,k}$$
$$\forall s(:STy_s = 'TF'), t(:t = 1..NT), k \tag{8}$$

$$\sum_{p(:SP_{s,p}=1 \& p \in P^F)}^{P} yPSl_{s,p,t,k} = yOp_{s,t,k}$$
$$\forall s(:STy_s = 'TF'), t(:t = 1..NT), k \tag{9}$$

**Mass Balance Constraints for raw water inventory**

The equations (10-11) represent mass balance on raw water inventory, accounting the total consumption of raw water received at the treatment plant for production of treated water products, as per throughput capacity and recovery yield of the treatment process. The reservoir capacity constraints (raw water and treated water inventory) limit the raw water supply from source to TF.

$$xQ_{s,i,p,t,k} = xQ_{s,i,p,t-1,k} + \sum_{s'(:SS_{s',s}=1 \text{ and } t-T_{s',s,p}^{RWTransit}>0)}^{S} xSSupl_{s',s,p,,k,t-T_{s',s,p}^{RWTransit}}$$
$$- \sum_{p'(:p' \in P^F)}^{P} yPsl_{s,p',t,k} * STpt_s * \frac{1}{\beta_{s,p'}}$$
$$\forall s(:STy_s = 'TF'), i(:i = 'RWI'), p(:p \in P^{RW} \text{ and } SIP_{s,i,p} = 1), t(:t = 2..NT), k \tag{10}$$

$$xQ_{s,i,p,t,k} = Q_{s,i,p,k}^{ini}$$
$$\forall s(:STy_s = 'TF'), i(:i = 'RWI'), p(:p \in P^{RW} \text{ and } SIP_{s,i,p} = 1), t(:t = 1), k \tag{11}$$

**Mass balance Constraints for treated water inventory**

The equations (12-13) represent mass balance on treated water product inventory in the treatment facility. Furthermore, equations (14-15) represent minimum and maximum capacity constraint for raw and treated water reservoirs respectively. The equations (16- 17) denotes buffer capacity and target capacity constraints in raw water and treated water reservoirs respectively for smooth operation of treatment plant and treated water supply operations.

$$xQ_{s,i,p,t,k} = xQ_{s,i,p,t-1,k} + yPSl_{s,p,t,k} * STpt_s - \sum_{c}^{C} xDeCon_{s,c,p,t}$$
$$\forall s(:STy_s = 'TF'), i(:i = 'TWI'), p(:p \in P^F \text{ and } SIP_{s,i,p} = 1), t(:t = 2..NT), k \quad (12)$$

$$xQ_{s,i,p,t,k} = Q^{ini}_{s,i,p,k}$$
$$\forall s(:STy_s = 'TF'), i(:i = 'TWI'), p(:p \in P^F \text{ and } SIP_{s,i,p} = 1), t(:t = 1), k \quad (13)$$

$$xQ_{s,i,p,t,k} \geq ICap^{min}_{s,i,p}$$
$$\forall s(:STy_s = 'TF'), i, p(:p \in P \text{ and } SIP_{s,i,p} = 1), t(:t = 1..NT), k \quad (14)$$

$$xQ_{s,i,p,t,k} \leq ICap^{max}_{s,i,p}$$
$$\forall s(:STy_s = 'TF'), i, p(:p \in P \text{ and } SIP_{s,i,p} = 1), t(:t = 1..NT), k \quad (15)$$

$$xQ_{s,i,p,t,k} \geq ICap^{buffer}_{s,i,p} - xBCV_{s,i,p,t,k}$$
$$\forall s(:STy_s = 'TF'), i(:i = 'RWI'), p(:p \in P^{RW} \text{ and } SIP_{s,i,p} = 1), t(:t = 1..NT), k \quad (16)$$

$$xQ_{s,i,p,t,k} - xTV^{+}_{s,i,p,t,k} + xTV^{-}_{s,i,p,t,k} = ICap^{Target}_{s,i,p,t}$$
$$\forall s(:STy_s = 'TF'), i(:i = 'TWI'), p(:p \in P^F \text{ and } SIP_{s,i,p} = 1), t(:ICap^{Target}_{s,i,p,t} > 0), k \quad (17)$$

**Ground water supply constraint**

The equation (18) represents constraint of maximum allowable extraction limit for ground water sources, usually specified by Central Ground Water Board, preventing the exploitation of the groundwater sources beyond replenishable limits.

$$xSSupl_{s,s',p,t,k} \leq SMax_s$$
$$\forall s(:SMax_s > 0), s'(:SS_{s,s'} = 1), p(:p \in P^{RW}), t(:t = 1..NT), k \quad (18)$$

**Consumer demand fulfillment constraint**

The equation (19) represents mass balance on consumer demand fulfillment from the total water supplied from all the sources.

$$\sum_{\substack{s(:SP_{s,p}=1 \text{and} SC_{s,c}=1 \\ \text{and} t-T^{Transit}_{s,c,p}>0)}}^{S} xDeCon_{s,c,p,t-T^{Transit}_{s,c,p}} + \Delta dem^{-}_{c,p,t,k} - \Delta dem^{+}_{c,p,t,k} = De_{c,p,t,k}$$

$$\forall c, p(: p \in P^F), t(: De_{c.p.t.k} > 0), k \tag{19}$$

**Tanker inventory**

The equations (20-21) represent the constraint that total distribution requirements should be equal to the total available distribution capacity of tankers in the demand horizon. Further, equation (22) represents overall time capacity balance constraint[32]. The equations (23 – 26) denotes the calculation of tanker capacity used in supplying raw water from the source to treatment plants and treated water to consumers in a particular time period. And, equation (27) represents hourly distribution capacity balance (T=1h being smallest time period unit in the formulation). The equation (28) balances tanker inventory along with any extra tanker hiring requirements at the start of planning horizon.

$$\sum_{t=1}^{NT} xDeCon_{s,c,p,t} - \sum_{v(:CPV_{c,p,v}=1)}^{V} xCDistb_{s,c,p,v,k} = 0$$
$$\forall c, p(: p \in P^F), s(: SC_{s,c}=1 \& SP_{s,p}=1), k \tag{20}$$

$$\sum_{t=1}^{NT} xSSupl_{s,s',p,t,k} - \sum_{v(:SSPV_{s,s'p,v}=1)}^{V} xVSSupl_{s,s',p,v,k} = 0$$
$$\forall s, s'(: SS_{s,s'}=1), p(: p \in P^{RW}), k \tag{21}$$

$$\sum_{c(:SP_{s,p}=1 \text{and} RS_{r,s}=1 \text{and} CPV_{c,p,v}=1)}^{C} \left( \frac{2 * T^{Travel}_{s,c,p,v} + T^{Prep}_{s,v} + T^{Disf}_{s,v} + T^{Distb}_{c}}{VQ_v} \right) * xCDistb_{s,c,p,v,k}$$
$$+ \sum_{s(:RS_{r,s}=1 \text{and} SSPV_{s,s'p,v}=1)}^{S} \left( \frac{2 * T^{Travel}_{s,c,p,v} + T^{Prep}_{s,v}}{VQ_v} \right) * xVSSupl_{s,s',p,v,k}$$
$$\leq NT * \left( VA_{r,v,p} + \frac{xVExQ_{r,v,p}}{VQ_v} \right)$$
$$\forall r, v, p(: RVP_{r,v,p}=1), k \tag{22}$$

$$\sum_{v(:CPV_{c,p,v}=1)}^{V} xPDl_{s,c,p,v,t,k} = xDeCon_{s,c,p,t}$$
$$\forall s, c(:SC_{s,c} = 1), p(:SP_{s,p} = 1 \& p \in P^F), t(:t = 1..NT), k \quad (23)$$

$$\sum_{t=1}^{NT} xPDl_{s,c,p,v,t,k} = xCDistb_{s,c,p,v,k}$$
$$\forall s, c(:SC_{s,c} = 1), p(:SP_{s,p} = 1 \text{ and } p \in P^F), v(:CPV_{c,p,v} = 1), k \quad (24)$$

$$\sum_{v(:SSPV_{s,s',p,v}=1)}^{V} xRw_{s,s',p,v,t,k} = xSSupl_{s,s',p,t,k}$$
$$\forall s, s'(:SS_{s,s'} = 1), p(:SSPV_{s,s',p,v} = 1), t(:t = 1..NT), k \quad (25)$$

$$\sum_{t=1}^{NT} xRw_{s,s',p,v,t,k} = xVSSupl_{s,s',p,v,k}$$
$$\forall s, s'(:SS_{s,s'} = 1), p(:SSPV_{s,s',p,v} = 1), v, k \quad (26)$$

$$xVQ_{r,v,p,t,k} = xVQ_{r,v,p,t-1,k} - \sum_{c}^{C} \sum_{s(:SP_{s,p}=1 \& RS_{r,s}=1 \& SC_{s,c}=1)}^{S} xPDl_{s,c,p,v,t,k}$$
$$- \sum_{s(:RS_{r,s}=1)}^{S} \sum_{s'(:SSPV_{s,s',p,v}=1)}^{S} xRw_{s,s',p,v,t,k}$$
$$+ \sum_{c}^{C} \sum_{s(:SP_{s,p}=1 \& RS_{r,s}=1 \& SC_{s,c}=1)}^{S} \sum_{\substack{t-1(:TE_{t-1}>(TE_{t'}+\frac{2*T_{s,s',p,v}^{Travel}+T_{s,v}^{Prep}+T_{s,v}^{Disf}+T_{c}^{Distb}}{24})) \\ t'(:TS_{t-1} \le (TS_{t'}+\frac{2*T_{s,s',p,v}^{Travel}+T_{s,v}^{Prep}+T_{s,v}^{Disf}+T_{c}^{Distb}}{24}))}} xPDl_{s,c,p,v,t',k}$$
$$+ \sum_{s(:RS_{r,s}=1)}^{S} \sum_{s'(:SSPV_{s,s',p,v}=1)}^{S} \sum_{\substack{t-1(:TE_{t-1}>(TE_{t'}+\frac{2*T_{s,s',p,v}^{Travel}+T_{s,v}^{Prep}}{24})) \\ t'(:TS_{t-1} \le (TS_{t'}+\frac{2*T_{s,s',p,v}^{Travel}+T_{s,v}^{Prep}}{24}))}} xRw_{s,s',p,v,t',k}$$
$$\forall r, v, p, t(:t = 2..NT), k \quad (27)$$

$$xVQ_{r,v,p}^{ini} = VQ_v * VA_{r,v,p} + xVExQ_{r,v,p}$$
$$\forall r, v, p(:RVP_{r,v,p} = 1) \quad (28)$$

**Objective function**

The equation (29) represents minimization of total expected operating cost as the objective function of the optimization problem.

$$xobj = \sum_{s}^{S}\sum_{c}^{C}\sum_{p(:SP_{s,p}=1)}^{P}\sum_{v}^{V} TrCost_{s,c,p,v}^{Distb} * xCDistb_{s,c,p,v}$$

$$+ \sum_{s}^{S}\sum_{s'(:SS_{s,s'}=1)}^{S}\sum_{p}^{P}\sum_{v(:SSPV_{s,s',p,v}=1)}^{V} TrCost_{s,s',p,v}^{RWsupply} * xVSSupl_{s,s',p,v}$$

$$+ \sum_{s}^{S}\sum_{i}^{I}\sum_{p}^{P} TVCost_{s,i,p} * \sum_{t=1}^{NT} xTV_{s,i,p,t} + \sum_{s}^{S}\sum_{i}^{I}\sum_{p}^{P} BCVCost_{s,i,p} * \sum_{t=1}^{NT} xBCV_{s,i,p,t}$$

$$+ \sum_{r}^{R}\sum_{v}^{V}\sum_{p(:RVP_{r,v,p}=1)}^{P} VExCost_{v,p} * \frac{xVExQ_{r,v,p}}{VQ_v}$$

$$+ \sum_{k}^{K} \pi_k \sum_{c}^{C}\sum_{p \in P^F}^{P}\sum_{t=1}^{NT}(q_p^+ \Delta dem_{c,p,t,k}^+ + q_p^- \Delta dem_{c,p,t,k}^-)$$

(29)

**3.3 Uncertainty Modelling:**

*3.3.1 Demand*: As discussed earlier in the section 2.2, the most widely used scenario tree method for discretized representation of uncertain parameters in stochastic optimization problems blows the size of the problem in an exponential manner with total number of scenarios (Tong et al., 2013). In such case, the problem becomes computationally intractable for solving in considerable time limit with a huge number of scenarios. On the flip side, a very small set of scenarios lacks full representation of uncertainty space. Apart from this, a major difficulty in solving scenario tree based stochastic optimization problem lies in computing the expectation value in the objective function. Therefore, we use the Monte Carlo simulation technique to create a finite reasonable number of scenario set based on the probability distribution function and coupled with (SAA) sample average approximation approach to approximate the expected total cost (using law of large numbers for *iid* samples). Hence, the demand uncertainty is represented in the formulation through *N* randomly generated scenarios (based on the normal distribution known from water consumption behaviour of customers) using MCS. These *N* discrete scenarios are represented by the index *k* and modelled in the constraints as Equation (30). Normal distribution for water demands is used herein [35] and the average value of the predicted demands ($De_{c,p,t,k}$) in Equation (30) are taken from the historical time series water consumption database. The discretized demand uncertainty parameter ξ_(*c,p,t,k*) in Equation (30) are taken as random seed for each scenario from standard normal probability distribution truncated to ± 20 % around the mean value.

$$xDe_{c,p,t,k} = De_{c,p,t}(1 + \xi_{c,p,t,k}) \quad \forall\, c, p, t, k \tag{30}$$

*3.3.2 Travel Time*: The uncertainty in travel time of tankers arising due to traffic congestions, loading operations at water source site or any unforeseen event related to tanker service from sources to consumers and treatment facilities is represented by two levels (nominal and high) with their assigned probabilities. Thus, corresponding to every MCS demand scenario in previous case (demand uncertainty), there are two nodes denoting travel time uncertainty realization in terms of nominal and high values, as shown in Fig.2. This gives a total of *2N* scenarios for this case of integrated demand and travel time uncertainty in the supply system.

## 4. Case Study

The efficacy of the developed stochastic optimization formulation in addressing the complexities and peculiarities of the tanker-based water distribution system operations is demonstrated herein through a representative case study.

**Input Data**: The proposed case study includes 4 water sources and 2 treatment facilities spanned across 3 regions of water supply in the city, as per the description in Table 1 (GW represents ground water source and FW represents fresh water source). The scheduling horizon considered is 5 days and is discretized on an hourly basis. All the input data to the model for this application case study is qualitatively corroborated with insights from literature papers, news-paper articles and survey reports on the tanker water systems (cited in introduction section) and understandings from discussions with a commercial enterprise - *Just Paani Water Supply Solutions*, to make it pertinent and representative of real-life problem and scale. Three type of water quality states are considered in this case study, namely, raw water (RW), domestic purpose water (DPW), and ultra-pure drinking water (UPDW). Furthermore, three type of consumer groups are considered, namely, house-hold use consumers (HHC), commercial scale consumers (CC) and hospital (institutional) consumers (HC), each demanding specific water quality product requirement as shown in Table 1. The information regarding advanced water treatment facilities throughput, recovery, reservoir capacities and all other related data is provided in the Supplementary information. The distance between consumer group locations, water sources and treatment facilities to calculate consumer-source and source- treatment facility suitability for tanker supply is also provided in supplementary information. In the proposed case study, we consider two type of tanker capacities, namely, 6KL and 10KL with different material trucks for supply of raw and treated water. The total initial availability of each tanker type is provided region wise as an input to the model. The full details of region wise availability along with speed, cost per kilo-meter and extra tanker hiring cost of each tanker type can be found in Supplementary information.

Table 1: Case Study problem specifics

| Region ID | Source/ TF ID | Consumer Type ID and (Water Product Demand) |
|---|---|---|
| R1 | GW1, TF1 | HHC1 (DPW & UPDW), CC1(DPW) |
| R2 | FW2, TF2 | CC2 (DPW), HC2 (UPDW) |
| R3 | FW3, GW3 | HHC3 (DPW, UPDW), HC3 (UPDW) |

**Scenario Generation**: The effect of demand and travel time uncertainty needs due consideration in the operational scheduling to develop a reliable tanker water supply system plan for efficient service. In this case study, we first consider only water demands of both water products (DPW and UPDW) for each consumer type as an uncertain parameter. Herein, based on the historical water consumption pattern, the demand uncertainty parameter ($\xi$) in Eq.(30) is estimated from a standard normal distribution in the range of ± 20% from the mean value. Further, in second case, where we consider both water demands and tanker travel time from source to consumers and treatment facilities as an uncertain parameter. Considering the uncertain delays in tanker travel durations owing to traffic congestions and poor road infrastructure in many urban (peri-urban) areas, the uncertainty in travel time is considered in two scenarios, namely NOMINAL and HIGH (being 30% high than nominal value of travel time). The probability of each travel time uncertainty scenario in this case is assumed to be 50%. As a result, corresponding to each demand uncertainty scenario there are two further scenarios of travel time uncertainty, resulting in the multiplied probability (independent events) in each *2N* scenarios.

Furthermore, the main objective of this framework is to supply good quality water to the consumers as much as possible and with that intent we assume a high penalty for demand shortfall in the case study presented. Secondly, based on discussions with a water tanker supply company named *Justpani*, we also want to capture the loss in business due to shortfall in demand supply. Keeping both the above objectives in mind, we have set the shortfall penalty as 100 times of the distribution cost. Furthermore, field data seems to suggest that the ratio of demand surplus to shortfall penalty costs can be assumed to be 1:10 for DPW (domestic purpose water) product supply (and 1:50 for UPDW product supply to hospital consumers) to minimize the instances of demand shortfall in all uncertainty scenarios. Moreover, we have performed a sensitivity analysis around the assumed surplus/shortfall cost penalty factors in the work. The figure below shows the sensitivity analysis plot (% age of scenarios in demand shortfall v/s surplus cases) where it can be observed that the assumed values (300/3000) of parameters results in 88 % surplus and only 9% demand shortage cases out of all possible uncertainty scenarios and, further increase in the value does not result in any significant impact on the percentage ratio of surplus and shortfall demand scenarios. Therefore, actual cost parameters have been set to be 300/3000. Similarly, cost penalty factors were chosen from sensitivity as 100/5000 for UPDW product to hospital consumers (figure is shown only for DPW product to prevent repetition and maintain the brevity of the manuscript).

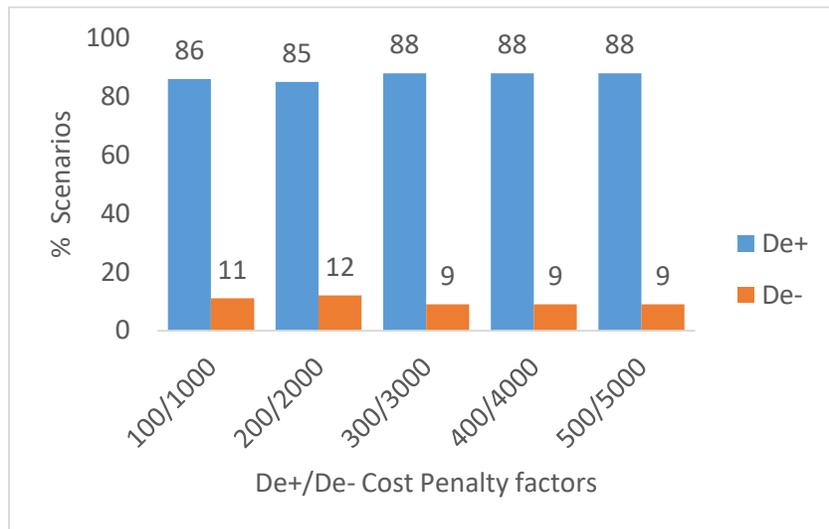

Fig 3. Sensitivity analysis for demand surplus and shortfall cost penalty factors in objective function

**Results and Discussion**: The optimization problem is programmed for solution in the FICO Xpress Optimization suit [36], using the "mmxprs" module 2.8.1, on a HP computer with an intel core i7, 16GB RAM PC. The Monte Carlo simulations for generating random demand scenarios were performed on MATLAB R2021. The case study discussed above are optimized for maximum of 36000s (10h) of CPU time or 1% optimality gap.

**4.1: Results: Only Demand Uncertainty Case**

Firstly, in order to determine the number of scenarios that can reasonably represent the demand uncertainty as well as computationally tractable, a convergence analysis is performed by increasing Monte Carlo sample scenarios from 5 to 60 in subsequent runs. The objective function value corresponding to each increment scenario of 5 samples is shown in Fig.4. As can be observed in Fig.3. the objective function value on y-axis starts stabilizing in the fluctuation band of 0.1 % when the sample scenarios is more than 45 and converges to 4.36 e+06 at $N$=50. Therefore, 50 scenarios of uncertain demand are considered in this case study. The equivalent deterministic model of this two-stage stochastic problem corresponding to $N$=50 scenarios have 24300 binary variables, 85144 continuous variables, 86161 constraints and took 434 sec to solve with 0.97% optimality gap.

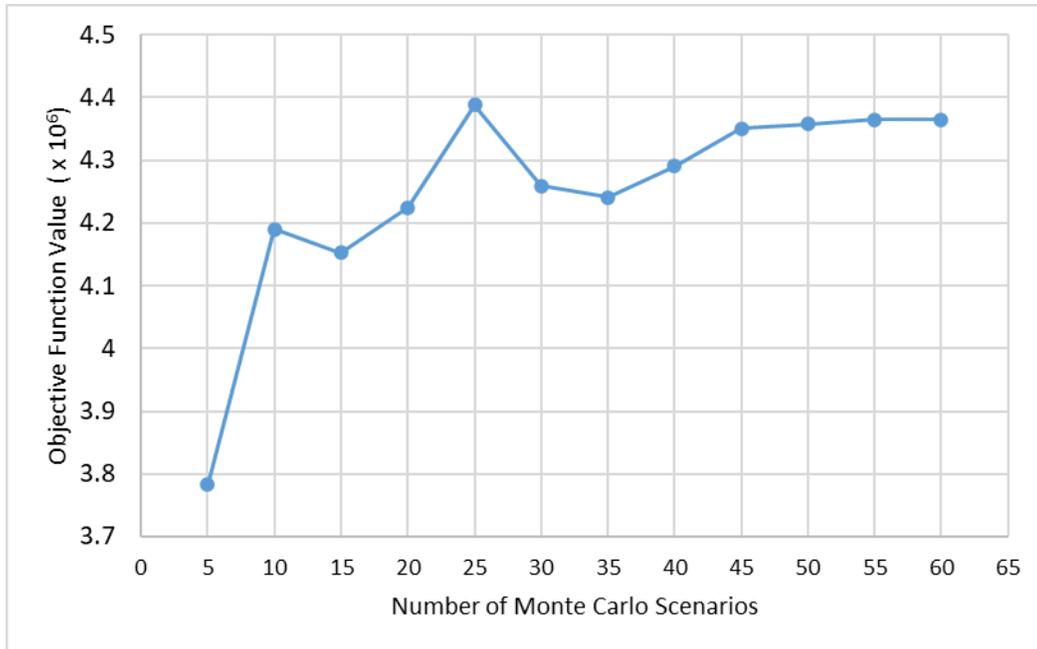

Fig 4. MCS convergence analysis

Herein, the optimal solution of two-stage stochastic problem with recourse approach is obtained by solving its equivalent deterministic MILP model. This solution value is abbreviated by (TSR). Secondly, if we could accurately estimate the actual water demands of a scenario *k* to be occurred in future (deterministic information of scenario *k*), the optimal solution corresponding to minimum operating cost can be obtained using purely deterministic model (DT). Furthermore, when the uncertain demand $D_{c,p,t,k}$ ($\xi$) is replaced by the expected mean value ($\bar{D}_{c,p,t}$), the corresponding deterministic model solution gives a first stage feasible solution for the two-stage stochastic problem. Fixing the first stage variables to this solution value and solving for recourse model for second stage decisions corresponding to all scenarios $k \in K$, a local optimum solution of the two-stage stochastic programming model can be found. This solution is called as expected demand value solution (EV).

Subsequently, a comparison of optimal results of total operating cost, total transportation cost in supplying water from various sources to WTPs and to customers via tankers, penalty costs for demand shortfall at consumer end, and extra tankers hired over the initial availability, in each model, namely, (i) deterministic (DT), (ii) expected demand value (EV) and, (iii) two-stage stochastic problem with recourse (TSR) is shown in Fig.5. The comparison shows that demand shortfall reduces significantly with TSR solution at a much lower operating cost than EV solution to prepare for demand uncertainty. This is because extra tanker hiring is a first stage decision and needs to be decided before we implement the schedule. However, EV solution has first stage decisions fixed (from average demand in DT model), it calculates the amount of demand which would be in shortfall/surplus category once the actual scenario surfaces on implementation of average demand solution. Hence, results of EV model solution in Fig.4 shows high shortfall cost as the tanker hiring is a fixed first stage decision (in this case study, no tankers were required for completely fulfilling the average demands in first stage).

Consequently, two measures of evaluating the benefits of performing stochastic optimization are calculated, namely (i) Value of Stochastic Solution (VSS), and (ii) Expected value of Perfect Information (EVPI) as per Eq. 6 and Eq. 7 respectively. The results of VSS and EVPI are shown in Table 2.

VSS= (TSR-EV)/EV  (6)

$$\text{EVPI} = (\text{TSR} - \text{DT})/\text{TSR} \tag{7}$$

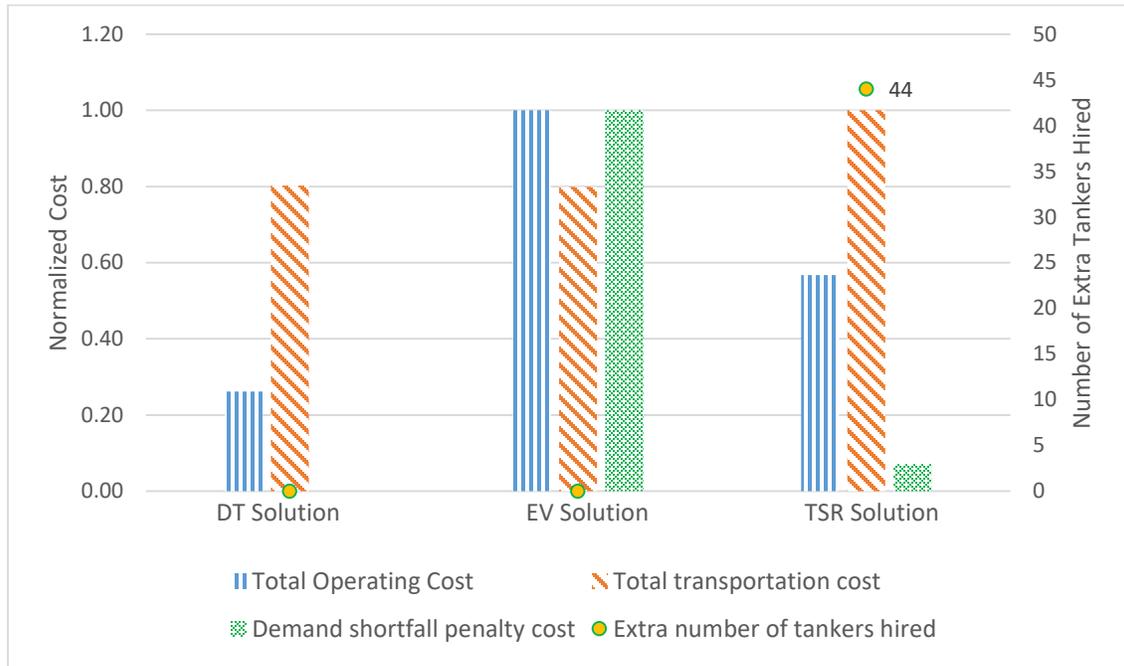

Fig.5. Comparison of different objective function cost components in optimal solution of three optimization models for demand uncertainty case

Table 2: VSS and EVPI Results for demand uncertainty case

| Objective Function Value | Deterministic Solution Model (DT) | Expected Value Solution Model (EV) | Two-Stage Stochastic Programming with Recourse Model (TSR) | Value of Stochastic Solution (VSS) | Expected Value of Perfect Information (EVPI) |
|---|---|---|---|---|---|
| Total Operating Cost | 2.01E+06 | 7.68E+06 | 4.36E+06 | -43 % | 54% |

As can be observed in Table 2, the total cost in stochastic with recourse approach solution is 43% lower than that of the expected value solution in this problem. Hence, in the face of real situations with uncertain demands, the VSS, i.e. difference between taking the average value as the solution as compared to performing stochastic analysis (propagating the uncertainties through the model and finding the impact on the objective function) for decision making results in significant cost savings. Therefore, the common strategy of using the expected value of the uncertain demand variable for optimization is found to be suboptimal. On the other hand, the deterministic (DT) model solution in Table 2 does not incorporate any uncertainties and assumes that perfect information is available. Thus, although the stochastic solution shows 54% increase in the operating cost compared with DT solution with perfect information, in real-life situations the demand and travel times are forecasted based on historical data; thereby getting perfect information about these parameters are unrealistic. Hence, stochastic solution actually provides an implementable solution while accounting the effect of uncertainties.

Furthermore, in Fig.6, the optimal manner of DPW demand satisfaction from stochastic solution based on all 50 scenarios for entire scheduling horizon is illustrated using box plot representation for all three regions. The central red line in the box plot in Fig.6 denotes the median demand value, and the two edges at bottom and top of the box denotes the first and third quartiles of all considered 50 scenarios demand data respectively. The black whiskers represent extreme points (outside the above range) in the data. The asterisk symbol indicates the expected value of demand fulfilment (optimal stochastic solution) considering all demand scenarios and corresponding probabilities. It can be seen in Fig.6. that TSR optimal solution can completely fulfil the DPW demands in more than 75% scenarios on all days of the scheduling horizon in each supply region.

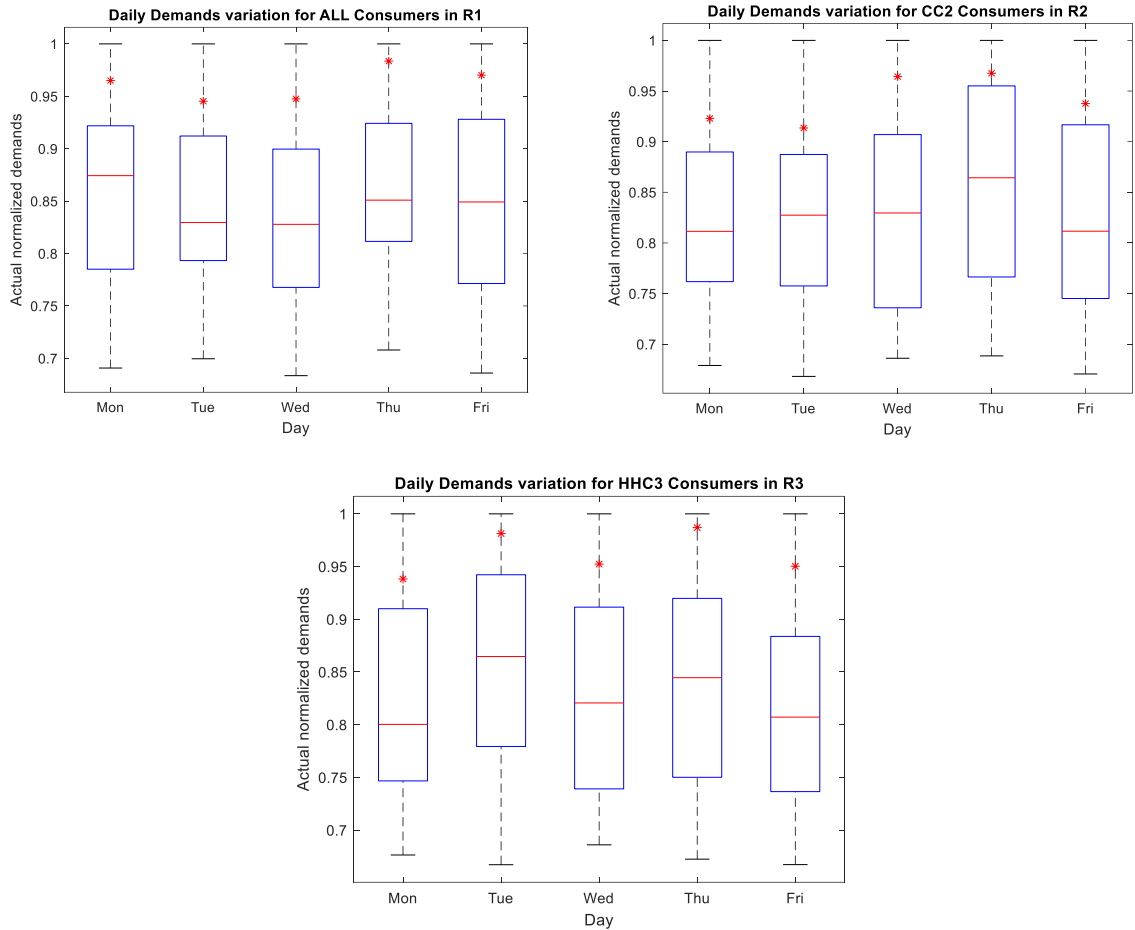

Fig. 6. DPW product demand satisfaction results based on 50 demand uncertainty scenarios in all three regions R1, R2, R3

Furthermore, Fig.7 shows mapping of sources and consumers for tanker water supply in deterministic and stochastic model solutions. The comparison between two cases clearly indicates shift of sources for minimum operating cost in the two cases to prepare for uncertain demands. Thus, this paper addresses the problem of maximizing fulfilment of demands and minimizing any demand shortfalls in water supply to consumers by tankers.

Fig.8 illustrates schedule of tanker movement for one selected day in the scheduling horizon (day 3) for both deterministic, and stochastic model for demand uncertainty. The green highlighted rows in Fig.8a and brown highlighted rows in Fig. 8b marks the difference in resultant optimal schedule for the two cases.

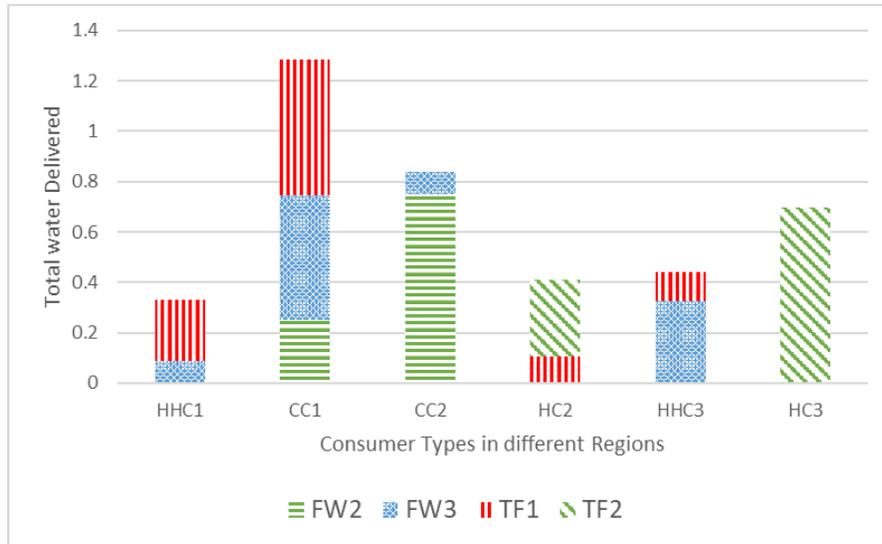

(a)

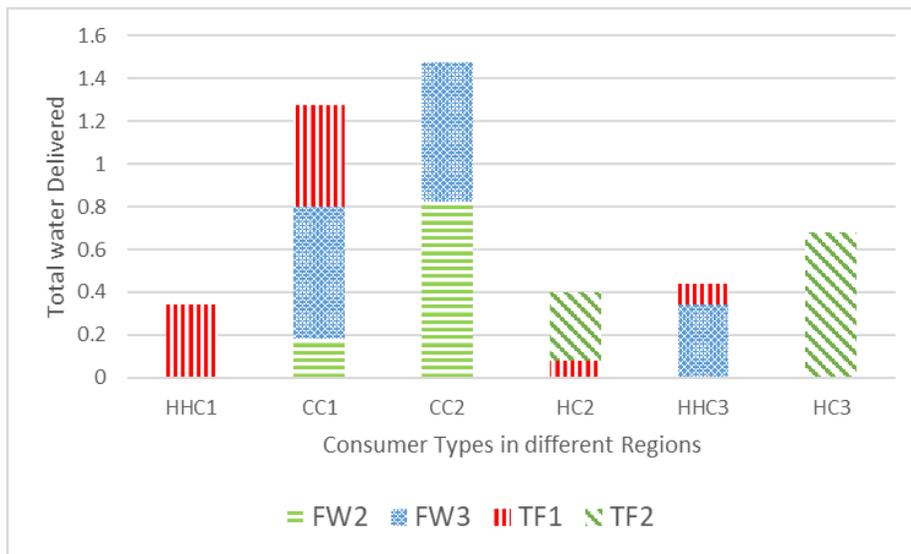

(b)

Fig.7: Source-Consumer mapping in optimal solution from two models (a) TSR model, (b) DT model

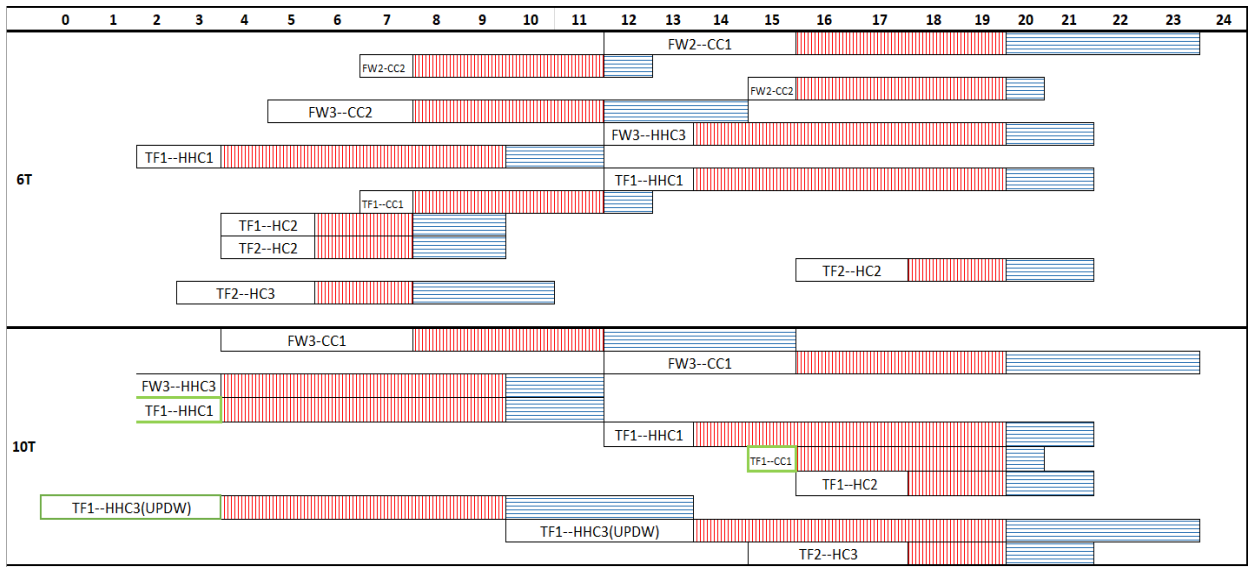

(a)

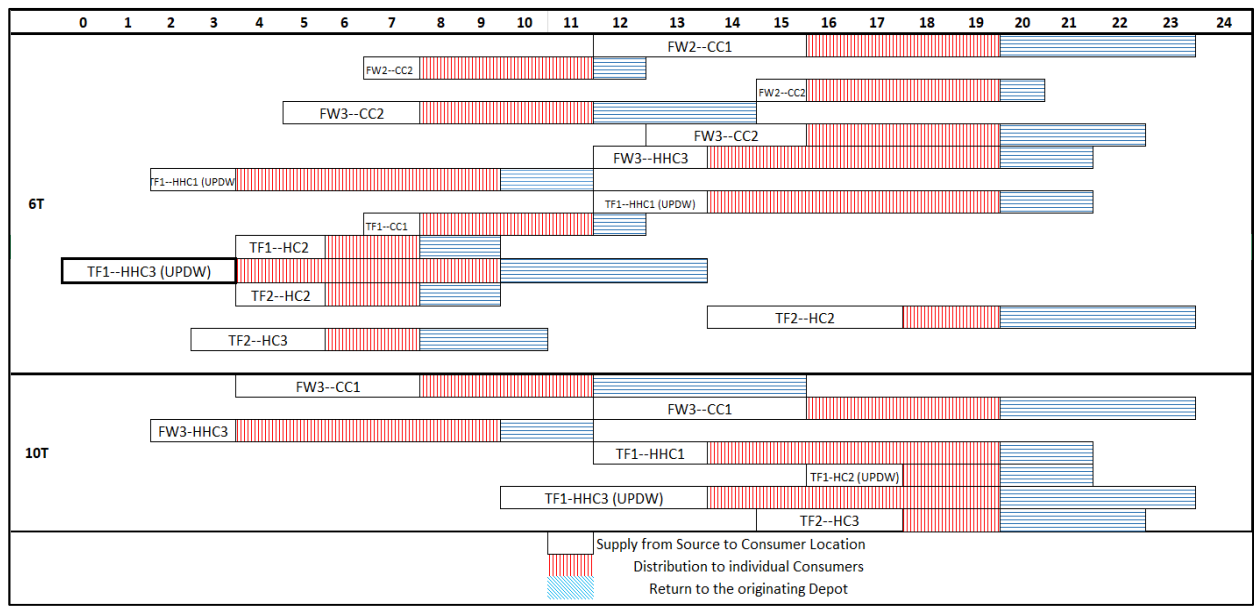

(b)

Fig.8. Schedule chart of water tanker movements on one selected day of scheduling horizon (day 3) in two cases (a) deterministic model (b) stochastic with demand uncertainty

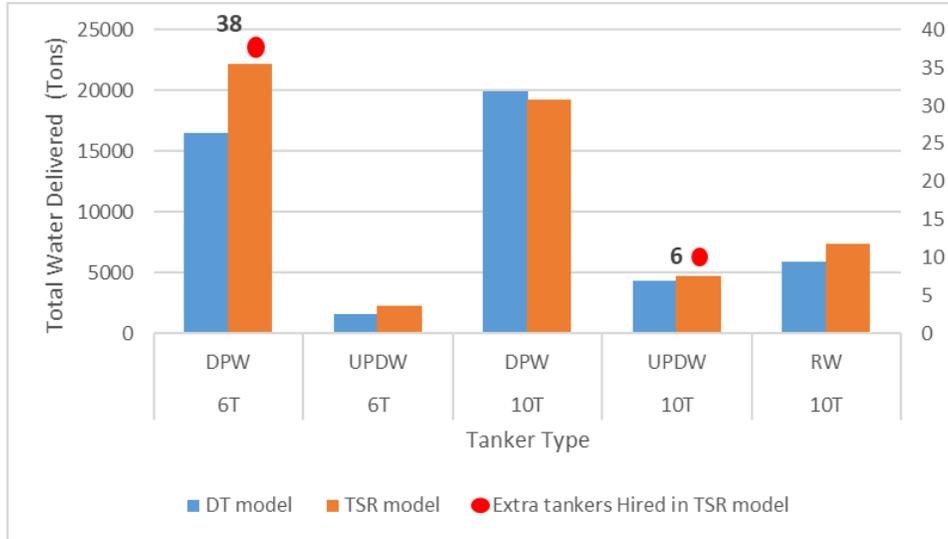

Fig.9. Total water delivered and extra tankers hired for each type for entire planning horizon

**4.2: Results: Integrated Demand and Travel Time Uncertainty Case**

Secondly, understanding of total water supply from different water sources in each case: (i) deterministic, (ii) stochastic model for only demand uncertainty and, (iii) stochastic model for demand plus travel time uncertainty is important. Fig.10 shows the contribution of each source (FW2, FW3, TF1 and TF2) in total demand fulfilment in the above mentioned three cases. It can be observed in Fig.9. that variability caused by integrated demand and travel time uncertainty is mostly catered by FW2 with its increased contribution (55%) in this case compared to 39% and 40% in deterministic and only demand uncertainty cases respectively. Furthermore, as travel time uncertainty becomes apparent in the integrated case and FW3 being far away from commercial consumer clusters (having highest DPW demands), dependence on it for water supply is reduced (24 %) in the optimal solution of this case compared to 31% and 33 % in deterministic and only demand uncertainty cases respectively. In addition, since the travel time uncertainty is prevalent in both routes, supply of water from GW sources to WTPs and from WTPs to consumers, as expected, contribution of TF1 is also decreased (12%) in integrated demand and time uncertainty case solution. Contribution of TF2 remain nominal in all three cases. Thus, the optimal contribution of each source in total water supply can be mapped with the proposed stochastic optimization formulation for each uncertainty case.

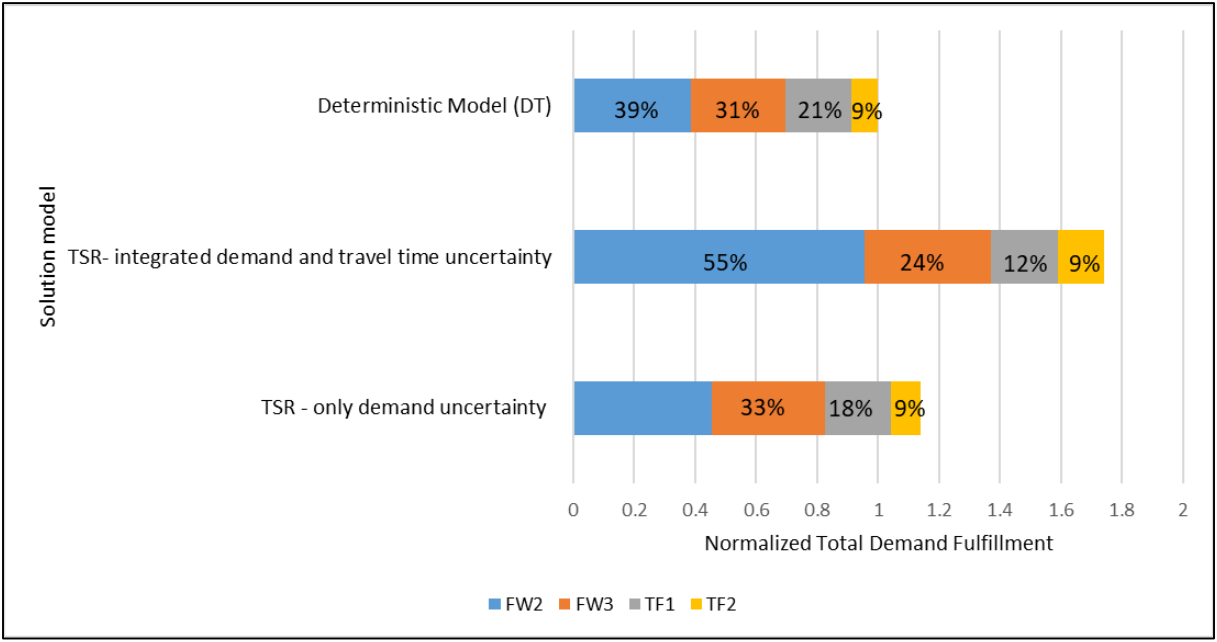

Fig.10. Contribution of each source in total water supply in three model results (i) Deterministic (DT), (ii) TSR for only demand uncertainty, and (iii) TSR for integrated demand and travel time uncertainty

Fig.11 shows objective function cost components compared for integrated demand and travel time uncertainty case and only demand uncertainty case. It can be seen in Fig.11 that total operating cost has increased significantly (12.32 x $10^6$) in integrated uncertainty case compared to (4.36 x $10^6$) in only demand uncertainty case optimal solution. As expected, this increase is corresponding to increased transportation cost, extra tanker hiring cost and demand shortage penalty cost in the integrated case. All these three components have increased cost due to incorporation of effect of time delays from travel time uncertainty, thus resulting in either requirement of a greater number of tankers or more number of tanker trips and increased penalty to minimize demand shortage instances in the total scheduling horizon.

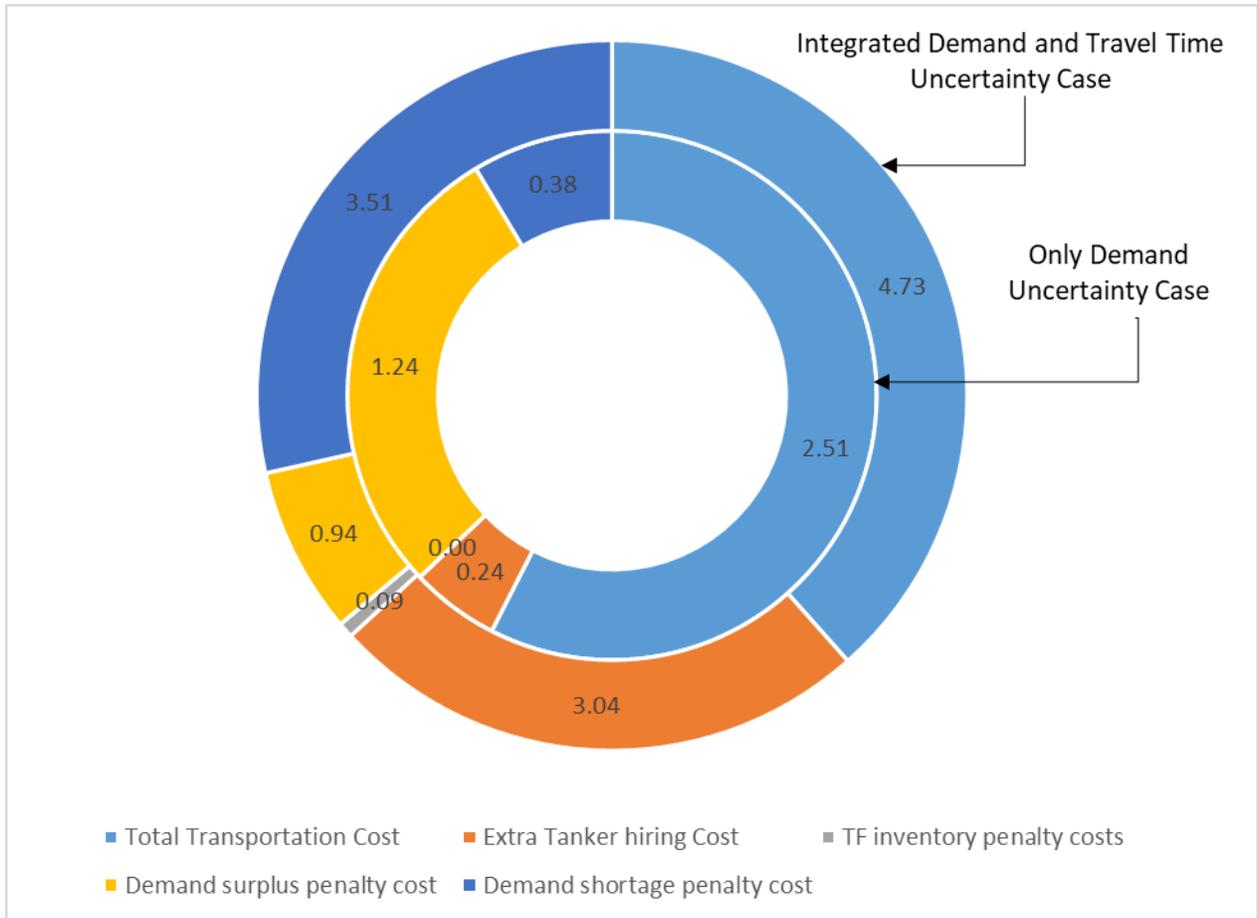

Fig.11. Comparison of Cost Components in only demand and integrated demand and travel time case solution objective function (Data value x 10^6)

5. **Conclusion**

Improving tanker water supply service management is increasingly considered as a matter of urban water supply planning. In this direction, understanding and preparing for impacts of demand and travel time uncertainty in tanker water supply system is crucial for efficient scheduling and management of day-to day operations. Furthermore, global water sustainability concerns demand not only the satisfaction of maximum water demands at the minimum operating cost but also the reduced wastage of water, limited extraction of ground water to allow it to replenish and supply of intended quality water through taker supply service.

To address these issues, a two-stage stochastic optimization framework with recourse approach is proposed in this paper, which aids in decision making related to water source selection, amount of water supply and tanker type from sources to consumers, tanker movement schedules and treatment facility operations. Uncertainty in consumer water demands is modelled using Monte Carlo simulation approach, while travel time uncertainty is accounted using scenario tree approach. The resultant hybrid scenario generation algorithm aids in reducing the two-stage stochastic programming problem into an equivalent deterministic optimization model with finite number of scenarios which makes the model computationally efficient. The consumer water demands are statistically represented by normal distribution based on the historical trend of water consumption behaviour in urban areas, whereas an empirical degree of belief is assigned to represent two scenarios of nominal and high case in travel time uncertainty, considering delays due to traffic congestions or any other unforeseen events while tanker filling and transportation.

The results indicated that stochastic optimization significantly improves the tanker water supply operations (in terms of minimizing cost, demand shortage penalty and reduced water wastages) compared to expected value solution of demand uncertainty. The proposed model also provides an optimal estimate of hiring extra tankers that need to fulfil the demands by exploiting the trade-off between cost of extra hired tankers, transportation cost and demand shortfall penalty. The developed framework is useful for informed decision making to assist in planning and scheduling of tanker-based water distribution system operations, considering uncertain inputs over the considered horizon. As a subject matter of future study, we also aim to address the uncertainty associated with seasonality of water source availability in the long-term planning formulation of tanker water distribution system.

## Acknowledgments

The authors gratefully acknowledge funding to the first author from the Ministry of Education, Government of India. Funding of the LOTUS project from Department of Science and Technology under the EU-India Water Co-operation program (Research and Innovation, Horizon 2020) is also gratefully acknowledged.

## Nomenclature

The following symbols and notation style is used in this paper:

**Acronyms:**

DPW    Domestic Purpose water

FW      Fresh water Source

GW      Ground water Source

RW      untreated raw water from ground water source

RWI     raw water inventory at treatment facility

TF      Treatment facility

TW      treated water in treatment facility

TWI     Treated water inventory at treatment plant

UPDW    Ultra-pure Drinking water

**Sets:**

| | |
|---|---|
| $C$ | set of consumers |
| $I$ | set of water inventory at treatment plants |
| $P$ | set of water product states |
| $P^{RW} \subseteq P$ | set of raw water states |
| $P^F \subseteq P$ | set of treated water state/final products |
| $R$ | set of regions for water supply in urban area |
| $S$ | set of water sources |
| $V$ | set of tanker vehicles |
| $WTP \subseteq S$ | set of water treatment plants |
| $K$ | set of demand scenarios |

**Indices:**

| | |
|---|---|
| $c \in C$ | consumer |
| $i \in I$ | water inventory in the treatment facility |
| $p \in P$ | water product state |

$r \in R$          region

$s, s' \in S$        water source

$v \in V$          tanker truck vehicle type

$t$            time period

$k \in K$          demand scenario

**Parameters**

$BCVCost_{s,i,p}$ penalty cost for violating the buffer capacity of product state $p$ in inventory $i$ at source $s$

$CPV_{c,p,v}$ assumes value equal to 1 if vehicle $v$ is compatible to deliver product of state $p$ to consumer $c$

$De^{\min}_{c,p,t}$ minimum aggregate demand from consumer group $c$ at time $t$ for water product of state $p$

$De^{\max}_{c,p,t}$ maximum aggregate demand from consumer group $c$ at time $t$ for water product of state $p$

$ICap^{\min}_{s,i,p}$ minimum capacity to be maintained at source $s$ in inventory $i$ for product of state $p$

$ICap^{\max}_{s,i,p}$ maximum capacity to be maintained at source $s$ in inventory $i$ for product of state $p$

$ICap^{buffer}_{s,i,p}$ buffer capacity limit at source s in inventory $i$ for product of state $p$

$ICap^{Target}_{s,i,p,t}$ target capacity at time $t$ for product of state $p$ in inventory $i$ at source $s$

$NT$ end time period of the planning horizon

$Op^{ini}_{s,t}$ captures operational state of treatment facility s at the start of the planning horizon

$Q^{ini}_{s,i,p}$ quantity of water product of state $p$ available in inventory $i$ of source $s$ at initial time period of planning horizon

$RS_{r,s}$ assumes value equal to 1 if source $s$ is suitable to transport water in region $r$

$RVP_{r,v,p}$ assumes value equal to 1 if tanker $v$ in region $r$ is compatible to supply product of state $p$

$SC_{s,c}$ assumes value equal to 1 if source $s$ is suitable to supply water to consumer $c$

$SIP_{s,i,p}$ assumes value equal to 1 for suitability of product of state $p$ with inventory $i$ at source $s$

$SMax_s$ Maximum groundwater extraction limit (KL/hour) from source $s$

$SS_{s,s'}$ Suitability of supplying raw water from source $s$ to treatment plant $s'$

$SSP_{s,s',p}$ assumes value 1 if source $s$ is suitable to supply product of state $p$ to treatment plant $s'$

$SSPV_{s,s',p,v}$ assumes value equal to 1 if vehicle type $v$ is compatible to supply product $p$ from source $s$ to treatment plant $s'$

$STpt_s$ throughput of water treatment plant $s$ to produce treated water

$STy_s$ indicates source type (FW/GW/TF)

$SP_{s,p}$ assumes value equal to 1 if source $s$ is suitable to supply water of state $p$

$T^{Disf}_{s,v}$ water disinfection time for a vehicle type $v$ from freshwater source $s$

$T^{Distb}_c$ Distribution time to consumers in a consumer group $c$

$T_s^{DT}$ treatment plant downtime period

$T_{s,v}^{\text{Pr}ep}$ preparation time for vehicle type $v$ at source $s$

$T_{s,c,p}^{Transit}$ transit time for transportation of water state $p$ from source $s$ to consumer $c$ (i.e. summation of preparation, disinfection and one way travel time)

$T_{s,s',p}^{RWTransit}$ transit time for transportation of water state $p$ from source $s$ to treatment plant s' (i.e. summation of preparation and one way travel time)

$T_{s,c,p,v}^{Travel}$ one way travel time required to deliver product state $p$ from source $s$ to consumer $c$ in vehicle type $v$

$T_s^{UT}$ treatment plant uptime period

$TE_t$ captures the start time of a time slot $t$

$TS_t$ captures the end time of a time slot $t$

$Tr\text{Cost}_{s,c,p,v}^{Distb}$ cost of transportation of unit quantity of water in tanker vehicle type $v$ for distributing treated water from source $s$ to consumer $c$

$Tr\text{Cost}_{s,s',p,v}^{RW\,supply}$ cost of transportation of unit quantity of water in tanker vehicle type $v$ for supplying raw water from source $s$ to treatment plant s'

$TV\text{Cost}_{s,i,p}$ penalty cost for violating the target capacity of fnal product state $p$ in inventory $i$ at source $s$

$VA_{r,v,p}$ number of tankers vehicles of type $v$ available in region $r$ suitable for transporting water product of state $p$

$VExCost_{v,p}$ penalty cost to purchase vehicle type $v$ for transporting product of state $p$

$VQ_v$ capacity of tanker vehicle type $v$

$\beta_{s,p}$ fraction of percentage recovery of permeate stream (treated water) from RO process at treatment plant $s$

$\xi_{c,p,t,k}$ random seed from standard normal distribution for demand uncertainty modelling

**Binary Decision Variables (prefix $y$ indicates all the binary variables)**

$yOp_{s,t}$ assumes value equal to 1 if treatment plant is operating at time period $t$ to produce treated water

$yPSl_{s,p,t}$ assumes value equal to 1 when product state $p$ is selected to be produced from treatment in treatment plant $s$'

**Continuous Decision Variables (prefix $x$ indicates all the continuous variables)**

$xBCV_{s,i,p,t}$ quantity of violation from buffer capacity limit of water product of state $p$ in inventory $i$ of source $s$ at time period $t$

$xCDistb_{s,c,p,v}$ quantity of water product of state $p$ supplied from source $s$ in tanker $v$ for distribution to consumer $c$

$xDeCon_{s,c,p,t}$ quantity of water product of state $p$ contributed by source $s$ at time $t$ to supply demand of consumer $c$

$xPDl_{s,c,p,v,t}$ quantity of water product of state $p$ delivered to consumer $c$ from source $s$ in tanker vehicle $v$ at time period $t$ $xQ_{s,i,p,t}$ quantity of water product of state *p* available in inventory *i* at time period *t* at source *s*

$xRw_{s,s',p,v,t}$ quantity of raw water supplied from source *s* to treatment plant *s'* in tanker vehicle type *v* at time period *t*

$xSDn_{s,t}$ assumes value equal to 1 if treatment plant *s'* is not operational at time period *t* to treat raw water

$xSSupl_{s,s',p,t}$ quantity of water product of state *p* supplied from source *s* to treatment plant *s'* at time period *t*

$xSUp_{s,t}$ assumes value equal to 1 if treatment plant *s'* is operational at time period *t* to treat raw water

$xTV_{s,i,p,t}$ quantity of violation from target limit of water product *p* in inventory *i* of source *s* at time period *t*

$xTV^{+}_{s,i,p,t}$ quantity of water product of state *p* in inventory *i* of source *s* at time period *t* which is positive violation from target value

$xTV^{-}_{s,i,p,t}$ quantity of water product of state *p* in inventory *i* of source *s* at time period *t* which is negative violation from target value

$xVExQ_{r,v,p}$ extra capacity of tanker vehicle type *v* required in region *r* for transporting product *p*

$xVSSupl_{s,s',p,v}$ quantity of water product of state *p* supplied in tanker vehicle type *v* from source *s* to treatment facility *s'*

$\Delta dem^{+}_{c,p,t,k}$ denote the surplus of the water product *p* for customer *c* at time *t* in scenario *k*

$\Delta dem^{-}_{c,p,t,k}$ denote the shortage of the water product *p* for customer *c* at time *t* in scenario *k*

For Table of Contents Only

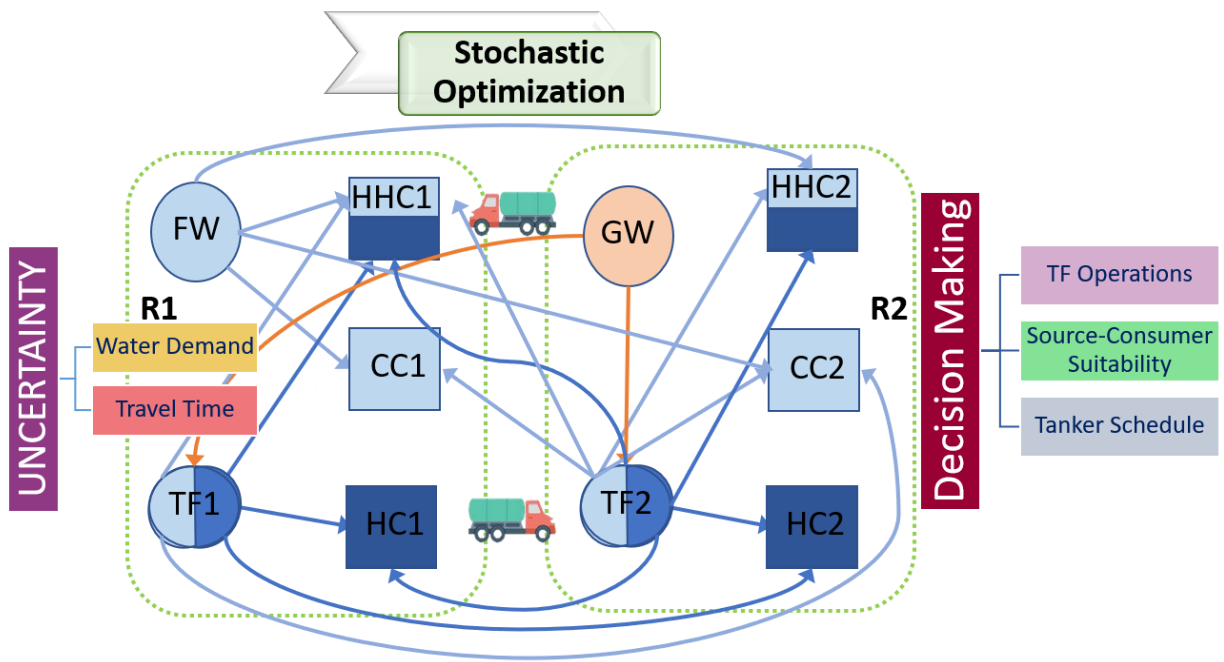